\documentclass[a4paper,draft,11pt]{article}
\usepackage{amsmath}
\newtheorem{theorem}{Theorem}[section]
\newtheorem{lemma}{Lemma}[section]
\newtheorem{corollary}{Corollary}[section]

\newtheorem{proposition}{Proposition}[section]

\newenvironment{remark}{\vspace{1ex}{\bf Remark.}\rm}{\vspace{1ex}}
\newenvironment{proof}{{\bf Proof.}}{\par\hspace{25em}\rule{1ex}{1ex}\par}

\newcommand{\ds}{\displaystyle}
\newcommand{\pt}{\partial}

\title{Totally geodesic submanifolds in the tangent bundle of a Riemannian 2-manifold.}
\author{Alexander YAMPOLSKY}
\begin{document}
\maketitle
\begin{abstract}
We give a full description of totally geodesic submanifolds in the tangent bundle of a Riemannian 2-manifold of
constant curvature and present a new class of a cylinder-type totally geodesic submanifolds in the general case.
   \\[2ex]
   {\it AMS subject class:} Primary 53B25; Secondary 53B20.
\end{abstract}


\section*{Introduction}
Let $(M^n,g)$ be a Riemannian manifold with metric $g$ and $TM^n$ its tangent bundle. S.~Sasaki \cite{Sk}
introduced on $TM^n$ a natural Riemannian metric $Tg$. With respect to this metric, all the fibers are totally
geodesic and intrinsically flat submanifolds. Probably M.-S. Liu \cite{Liu} was the first who noticed that the
base manifold embedded into $TM^n$ by the zero section is totally geodesic, as well. Soon afterwards, Sato K.
\cite{St} described geodesics (the totally geodesic submanifolds of dimension 1) in the tangent bundle over
space forms. The next step was made by P.Walczack \cite{Wch} who tried to find a non-zero section $\xi:M^n\to
TM^n$ such that the image $\xi(M^n)$ is a totally geodesic submanifold. He proved that if $\xi$ is of constant
length and $\xi(M^n)$ is totally geodesic, then $\xi$ is a parallel vector field. As a consequence, the base
manifold should be reducible. The irreducible case stays out of considerations up to now. A general conjecture
stated by A.Borisenko claims that, in irreducible case, the zero vector field is the unique one which generates
a totally geodesic submanifold $\xi(M^n)$ or, equivalently, the base manifold is the unique totally geodesic
submanifold of dimension $n$ in $TM^n$ transversal to fibers. A dimensional restriction is essential.
M.T.K.~Abbassi and the author \cite{Ab-Ym} treated the case of fiber transversal submanifolds in $TM^n$ of
dimension $l<n$ and have found some examples of totally geodesic submanifolds of this type. Earlier this problem
had been considered in \cite{Ym1}.

It is also worthwhile to mention that in the case of \textit{tangent sphere bundle } the situation is different.
Sasaki S. \cite{Sk1} described geodesics  in the tangent sphere bundle over space forms and Nagy P. \cite{Ng}
described geodesics in the tangent sphere bundle over symmetric spaces. The author has given a full description
of totally geodesic vector fields on 2-dimensional manifolds of constant curvature \cite{Ym2} and an example of
a totally geodesic unit vector field on positively/negatively curved manifolds of non-constant curvature
\cite{Ym3}. A full description of 2-manifolds which admit a totally geodesic unit vector field was given in
\cite{Ym4}.

In this paper we consider a more general problem concerning the description of all possible totally geodesic
submanifolds in the \textit{tangent bundle} of Riemannian 2-manifold with a sign-preserving curvature. For the
spaces of constant curvature this problem was posed by A.Borisenko in \cite{BY}.

 In Section \ref{2-dim} we prove the following theorems. \vspace{2ex}

 \textbf{Theorem 1.}\textit{ Let $M^2$ be Riemannian manifold of constant curvature
$K\ne0$. Suppose that $\tilde F^2\subset TM^2$ is a  totally geodesic submanifold. Then locally $\tilde F^2$ is
one of the following submanifolds:
\begin{itemize}
\item[(a)] a single fiber $T_qM^2$; \item[(b)] a cylinder-type surface based on a geodesic $\gamma$ in $M^2$
with elements generated by a parallel unit vector field along $\gamma$; \item[(c)] the base manifold embedded
into $TM^2$ by zero vector field.
\end{itemize}}
\vspace{2ex}

Remark that the item (b) of Theorem 1 is a consequence of more general result. \vspace{2ex}

\textbf{Theorem 2} \textit{Let $M^2$ be a Riemannian manifold of sign-preserving curvature. Suppose that $\tilde
F^2\subset TM^2$ is a totally geodesic submanifold having non-transversal intersection with the fibers. Then
locally  $\tilde F^2$ is a cylinder-type surface based on a geodesic $\gamma$ in $M^2$ with elements generated
by a parallel unit vector field along $\gamma$.}\vspace{2ex}

Moreover, a general Riemannian manifold $M^n$ admits this class of totally geodesic surfaces in $TM^n$ (see
Proposition \ref{Gener}).

In Section \ref{3-dim} we prove the following general result.\vspace{2ex}

\textbf{Theorem 3. }\textit{Let $M^2$ be a Riemannian manifold with sign-preserving curvature. Then $TM^2$ does
not  admit a totally geodesic 3-manifold even locally.}\vspace{2ex}

\textbf{Acknowledgement. } The author expresses his thanks to professor E. Boeckx (Leuven, Belgium) for useful
remarks in discussing the results.

\section{Necessary facts about the Sasaki metric}
Let $(M^n,g)$ be an $n$-dimensional Riemannian manifold  with metric $g$. Denote by $\big<\cdot\,,\cdot\big>$
the scalar product with respect to $g$. The {\it Sasaki metric} on $TM^n$ is defined by the following scalar
product: if $\tilde X,\tilde Y$ are tangent vector fields on $TM^n$,  then
\begin{equation}
\label{Eqn1}
  \big<\big<\tilde X,\tilde Y\big>\big>:=\big<\pi_* \tilde X, \pi_* \tilde Y\big>+\big<K \tilde X,K \tilde Y\big>,
\end{equation}
where $\pi_*:TTM^n \to TM^n $ is the differential of the projection $\pi:TM^n \to M^n $ and $K: TTM^n \to TM^n$
is the {\it connection map} \cite{Dmb}. The local representations for $\pi_*$ and $K$ are the following ones.

Let $(x^1,\dots ,x^n)$ be a local coordinate system on $M^n$. Denote by $\pt_i:=\partial/\partial x^i $ the
natural tangent coordinate frame. Then, at each point $q\in M^n$, any tangent vector $\xi$ can be decomposed as
$\xi=\xi^i \,\pt_i|_q$. The set of parameters $\{x^1,\dots ,x^n;\,\xi^1,\dots,\xi^n\}$ forms the natural induced
coordinate system in $TM^n$, i.e. for a point $Q=(q,\xi )\in TM^n$, with $q\in M^n, \ \ \xi \in T_qM^n$, we have
$q=(x^1,\dots ,x^n), \, \xi =\xi ^i\,\pt_i|_q$. The natural frame in $T_{Q}TM^n$ is formed by
$$
\tilde\pt_i:=\frac{\partial}{\partial x^i}|_Q, \quad \tilde\pt_{n+i}:=\frac{\partial}{\partial \xi ^i}|_Q
$$
and for any $\tilde X\in T_{Q}TM^n$ we have the decomposition
$$
\tilde X=\tilde X^i\tilde\pt_i+\tilde X^{n+i}\tilde\pt_{n+i}\ .
$$
 Now locally, the \textit{horizontal} and \textit{vertical} projections of $\tilde X$ are given by
   \begin{equation} \label{Eqn2}
   \begin{array}{l}
   \pi_* \tilde X|_Q= \tilde X^i\,\pt_i|_q, \\[1ex]
   K \tilde X|_Q= (\tilde X^{n+i}+ \Gamma^i_{jk}(q)\,\xi^j \tilde X^k)\, \pt_i|_q, \\[1ex]
   \end{array}
   \end{equation}
where $ \Gamma^i_{jk}$ are the Christoffel symbols of the metric $g$.

The inverse operations are called \textit{lifts }. If $  X=  X^i\,\pt_i$ is a vector field on $M^n$ then the
vector fields on $TM$ given by
   $$
   \begin{array}{l}
     X^h=  X^i\,\tilde\pt_i- \Gamma^i_{jk}\,\xi^j  X^k\,\tilde\pt_{n+i},\\[1ex]
     X^v=  X^i\,\tilde\pt_{n+i}
   \end{array}
   $$
are called the \textit{horizontal} and \textit{vertical} lifts of $X$ respectively. Remark that for any vector
field $  X$ on $M^n$ it holds
   \begin{equation}\label{Pr}
   \begin{array}{ll}
   \pi_* {  X}^h=  X,& K {  X}^h=0, \\[1ex]
   \pi_* {  X}^v=0,       & K {  X}^v=  X.
   \end{array}
   \end{equation}

There is a natural decomposition
$$
T_Q(TM^n)=\mathcal{H}_Q(TM^n)\oplus \mathcal{V}_Q(TM^n),
$$
where $\mathcal{H}_Q(TM^n)=\ker K$ is called the \textit{horizontal distribution} and $ \mathcal{V}_Q(TM^n)=\ker
\pi_*$ is called the \textit{vertical distribution} on $TM^n$. With respect to the Sasaki metric, these
distributions are mutually orthogonal. The vertical distribution is \textit{integrable} and the fibers are
precisely its integral submanifolds. The horizontal distribution is \textit{never integrable} except the case of
a flat base manifold.

{\it For any vector fields $  X,  Y$ on $M^n$, the covariant derivatives of various combinations of lifts to the
point $Q=(q,\xi) \in TM^n$ can be found by the formulas} \cite{Kow}

\begin{equation}\label{Kow}
\begin{array}{ll}
\tilde \nabla_{X^h}Y^h|_Q = (  \nabla_{ X}Y|_q)^h- \frac{1}{2}( R_q(X,Y)\xi)^v, \
&\tilde \nabla_{X^v}Y^h|_Q = \frac{1}{2}( R_q(\xi , X)Y)^h,\\[2ex]
\tilde \nabla_{X^h}Y^v|_Q = ( \nabla_{X}Y|_q)^v+ \frac{1}{2}(R_q(\xi ,Y)X)^h, \ & \tilde \nabla_{ X^v} Y^v|_Q=0,
\end{array}
\end{equation}
{\it where $  \nabla$ and $  R$ are the Levi-Civita connection and the curvature tensor of $M^n$ respectively}.

\begin{remark}
The formulas \eqref{Kow} are applicable to the \emph{lifts} of vector fields only. A formal application to a
general field on tangent bundle may lead to wrong result.  For example,
$$
\begin{array}{rl}
  \tilde\nabla_{X^v} (\xi^i (\partial_i)^h)
  =&  X^v(\xi^i)\, \partial_i^h
     + \xi^i \tilde\nabla_{X^v}\partial_i^h \\
  =&   X^i \partial_i^h + \xi^i \frac{1}{2} \big(R(\xi, X)\partial_i\big)^h
  = X^h + \frac{1}{2} \big(R(\xi, X)\xi \big)^h
\end{array}
$$
and we have an additional term in the formulas. We will use this rule in our calculations without special
comments.

\end{remark}

\section{Local description of 2-dimensional totally geodesic submanifolds in $TM^2$}\label{2-dim}

In this section we prove Theorem 1. The proof is given in a series of subsections. Namely, in subsection
\ref{Prel} we prove the item (a), in subsection \ref{VF} we prove the item (c) and finally, in subsection
\ref{Ruled} we prove Theorem 2 and therefore, the item (b) of Theorem 1.

\subsection{Preliminary considerations.}\label{Prel}
Let $\tilde F^2$ be a submanifold in $TM^2$. Let $(x^1,x^2;\xi^1,\xi^2)$ be a local chart on $TM^2$. Then
locally $\tilde F^2$ can be given by mapping $f$ of the form
$$
f: \left\{
\begin{array}{l}
x^1=x^1(u^1,u^2),\\[1ex]
x^2=x^2(u^1,u^2),
\end{array}
\right.
\quad
\begin{array}{l}
\xi^1=\xi^1(u^1,u^2),\\[1ex]
\xi^2=\xi^1(u^1,u^2),
\end{array}
$$
where $u^1,u^2 $ are the local parameters on $\tilde F^2$. The Jacobian matrix $f_*$ of the mapping $f$ is of
the form {\large
$$
f_*= \left(
\begin{array}{cc}
 \frac{\pt x^1}{\pt u^1} & \frac{\pt x^1}{\pt u^2}    \\[1ex]
 \frac{\pt x^2}{\pt u^1} & \frac{\pt x^2}{\pt u^2}    \\[1ex]
\frac{\pt \xi^1}{\pt u^1} & \frac{\pt \xi^1}{\pt u^2}   \\[1ex]
 \frac{\pt \xi^2}{\pt u^1} & \frac{\pt \xi^2}{\pt u^2} \\[1ex]
\end{array}
\right).
$$
} Since $rank \ f_*=2$, we have three \textit{geometrically different} possibilities to achieve the rank, namely
{\large
$$
\begin{array}{l}
(a)\quad \det
\left(%
\begin{array}{cc}
 \frac{\pt x^1}{\pt u^1} & \frac{\pt x^1}{\pt u^2} \\[1ex]
 \frac{\pt x^2}{\pt u^1} & \frac{\pt x^2}{\pt u^2} \\[1ex]
\end{array}
\right)\ne0; \qquad (b)\quad \det
 \left(%
\begin{array}{cc}
 \frac{\pt x^1}{\pt u^1} & \frac{\pt x^1}{\pt u^2}    \\[1ex]
\frac{\pt \xi^1}{\pt u^1} & \frac{\pt \xi^1}{\pt u^2}   \\[1ex]
\end{array}%
\right)\ne0;
\\[4ex]
(c)\quad \det
 \left(%
\begin{array}{cc}
 \frac{\pt \xi^1}{\pt u^1} & \frac{\pt \xi^1}{\pt u^2}    \\[1ex]
\frac{\pt \xi^2}{\pt u^1} & \frac{\pt \xi^2}{\pt u^2}   \\[1ex]
\end{array}%
\right)\ne0.
\end{array}
$$
} Without loss  of generality we can consider these possibilities in a way that (b) excludes (a), and  (c)
excludes (a) and (b) restricting the considerations to a smaller neighbourhood or even to an open and dense
subset.\vspace{2ex}

\textbf{Case (a).} In this case one can locally parameterize the submanifold under consideration as
$$
f: \left\{
\begin{array}{l}
x^1=u^1,\\[1ex]
x^2=u^1,
\end{array}
\right. \quad
\begin{array}{l}
\xi^1=\xi^1(u^1,u^2),\\[1ex]
\xi^2=\xi^2(u^1,u^2),
\end{array}
$$
and we can consider the submanifold $\tilde F^2$ as an image of the vector field $\xi(u^1,u^2)$ on the base
manifold. Denote $\tilde F^2$ in this case by $\xi(M^2)$. We analyze this case in subsection
\ref{VF}.\vspace{2ex}

\textbf{Case (b).} In this case one can parameterize the submanifold $F^2$ as
$$
f: \left\{
\begin{array}{l}
x^1=u^1,\\[1ex]
x^2=x^2(u^1,u^2),
\end{array}
\right. \quad
\begin{array}{l}
\xi^1=u^2,\\[1ex]
\xi^2=\xi^2(u^1,u^2).
\end{array}
$$
Taking into account that we exclude the case (a) in  considerations of the case (b), we should set
$$
\det \left(
\begin{array}{cc}
 \frac{\pt x^1}{\pt u^1} & \frac{\pt x^1}{\pt u^2} \\[1ex]
 \frac{\pt x^2}{\pt u^1} & \frac{\pt x^2}{\pt u^2} \\[1ex]
\end{array}
\right)= \det \left(
\begin{array}{cc}
 1 & 0 \\[1ex]
 \frac{\pt x^2}{\pt u^1} & \frac{\pt x^2}{\pt u^2} \\[1ex]
\end{array}
\right)=\frac{\pt x^2}{\pt u^2}=0.
$$
Therefore, $x^2(u^1,u^2)$ does not depend on $u^2$ and the local representation takes the form
$$
f: \left\{
\begin{array}{l}
x^1=u^1,\\[1ex]
x^2=x^2(u^1),
\end{array}
\right. \quad
\begin{array}{l}
\xi^1=u^2,\\[1ex]
\xi^2=\xi^2(u^1,u^2).\\[1ex]
\end{array}
$$
Remark that $\pi(\tilde F^2)= (u^1,x^2(u^1)$ is a regular curve on $M^2$. If we denote this projection by
$\gamma(s)$ parameterized by the arc-length parameter and set $u^2:=t$, the local parametrization of $\tilde
F^2$ takes the form
\begin{equation}\label{r_def}
\gamma(s): \left\{
\begin{array}{l}
x^1=x^1(s),\\[1ex]
x^2=x^2(s),
\end{array}
\right.
\qquad
 \xi(t,s): \left\{
\begin{array}{l}
\xi^1=t,\\[1ex]
\xi^2=\xi^2(t,s)
\end{array}
\right.
\end{equation}
We can interpret this kind of submanifolds in $TM^2$ as a one-parametric family of smooth vector fields over a
regular curve on the base manifold. We will refer to this kind of submanifolds as \textit{ruled submanifolds} in
$TM^2$ and analyze their totally geodesic property in subsection \ref{Ruled}.\vspace{2ex}

\textbf{Case (c).} It this case a local parametrization of $\tilde F^2$ can be given as
$$
f: \left\{
\begin{array}{l}
x^1=x^1(u^1,u^2),\\[1ex]
x^2=x^2(u^1,u^2),
\end{array}
\right. \quad
\begin{array}{l}
\xi^1=u^1,\\[1ex]
\xi^2=u^2.\\[1ex]
\end{array}
$$
Taking into account that we exclude the case (b) considering the case (c), we should suppose
$$
\det \left(
\begin{array}{cc}
 \frac{\pt x^1}{\pt u^1} & \frac{\pt x^1}{\pt u^2} \\[1ex]
 \frac{\pt \xi^1}{\pt u^1} & \frac{\pt \xi^1}{\pt u^2} \\[1ex]
\end{array}
\right)= \det \left(
\begin{array}{cc}
\frac{\pt x^1}{\pt u^1} & \frac{\pt x^1}{\pt u^2} \\[1ex]
 1 & 0 \\[1ex]
\end{array}
\right)=-\frac{\pt x^1}{\pt u^2}=0.
$$

$$
\det \left(
\begin{array}{cc}
 \frac{\pt x^1}{\pt u^1} & \frac{\pt x^1}{\pt u^2} \\[1ex]
 \frac{\pt \xi^2}{\pt u^1} & \frac{\pt \xi^2}{\pt u^2} \\[1ex]
\end{array}
\right)= \det \left(
\begin{array}{cc}
\frac{\pt x^1}{\pt u^1} & \frac{\pt x^1}{\pt u^2} \\[1ex]
 0 & 1 \\[1ex]
\end{array}
\right)=\frac{\pt x^1}{\pt u^1}=0.
$$
Thus, we conclude $x^1=const$. In the same way, we get $x^2=const$. Therefore, \textit{a submanifold of this
kind is nothing else but the fiber, which is evidently totally geodesic and there is nothing to prove.}

\subsection{Totally geodesic vector fields}\label{VF}

 In \cite{Ab-Ym} the author has found the conditions on a vector field to generate a totally geodesic
submanifold in the tangent bundle. Namely, let $\xi$ be a vector field on $M^n$. The submanifold $\xi(M^n)$ is
totally geodesic in $TM^n$ if and only if for any vector fields $X,Y$ on $M^n$ the following equation holds
\begin{equation}\label{Eq1}
r(X,Y)\xi+r(Y,X)\xi -\nabla_{h_\xi(X,Y)}\xi=0,
\end{equation}
where $r(X,Y)\xi=\nabla_X\nabla_Y\xi-\nabla_{\nabla_XY}\xi$ is  "half" the Riemannian curvature tensor and
$h_\xi(X,Y)=R(\xi,\nabla_X\xi)Y+R(\xi,\nabla_Y\xi )X$.

It is natural to rewrite this equations in terms of $\rho$ and $e_\xi$ where $e_\xi$ is a \textit{unit} vector
field and $\rho$ is the length function of $\xi$.

\begin{lemma}
Let $\xi=\rho\,e_\xi$ be a vector field on a Riemannian manifold $M^n$. Then $\xi(M^n)$ is totally geodesic in
$TM^n$ if and only if for any vector field $X$ the following equations hold
\begin{equation}\label{Eq2}
\left\{
\begin{array}{l}
 Hess_\rho(X,X)-\rho^2\Big(R(e_\xi,\nabla_Xe_\xi)X\Big)(\rho)-\rho\,|\nabla_Xe_\xi|^2=0,\\[2ex]
 \rho^3\nabla_{R(e_\xi,\nabla_Xe_\xi)X}e_\xi-2X(\rho)\nabla_Xe_\xi-\rho(r(X,X)e_\xi+
 |\nabla_Xe_\xi|^2\,e_\xi)=0,
\end{array}
\right.
\end{equation}
where $Hess_\rho(X,X)$ is the Hessian of the function $\rho$.
\end{lemma}
\begin{proof}
Indeed, the equation (\ref{Eq1}) is equivalent to
\begin{equation}\label{Eq3}
r(X,X)\xi=\nabla_{R(\xi,\nabla_X\xi)X}\xi,
\end{equation}
where $X$ is an arbitrary vector field. Setting $\xi=\rho\,e_\xi$, where $e_\xi$ is a unit vector field, we have
$$
\begin{array}{ll}
r(X,X)\xi=&\nabla_X\nabla_X(\rho\,e_\xi)-\nabla_{\nabla_XX}(\rho\,e_\xi)=\\[1ex]
&\nabla_X(X(\rho)e_\xi+\rho\,\nabla_Xe_\xi)-(\nabla_XX)(\rho)e_\xi-\rho\,\nabla_{\nabla_XX}e_\xi=\\[1ex]
&\Big(X(X(\rho))-(\nabla_XX)(\rho)\Big)e_\xi+2X(\rho)\nabla_Xe_\xi+\rho\,r(X,X)e_\xi
\end{array}
$$
and
$$
\nabla_{R(\xi,\nabla_X\xi)X}\xi=\rho^2\Big(R(e_\xi,\nabla_Xe_\xi)X\Big)(\rho)\, e_\xi
+\rho^3\nabla_{R(e_\xi,\nabla_Xe_\xi)X}e_\xi.
$$
If we remark  that $X(X(\rho))-(\nabla_XX)(\rho)\stackrel{def}{=}Hess_\rho(X,X)$ and for a unit vector field
$e_\xi$
$$
\big<r(X,X)e_\xi,e_\xi\big>=-|\nabla_Xe_\xi|^2,
$$
then we can easily decompose the equation (\ref{Eq3}) into components, parallel to and orthogonal to $e_\xi$,
which gives the equations (\ref{Eq2}).
\end{proof}

\begin{corollary}\label{Cor1}
Suppose that $M^n$ admits a  totally geodesic vector field $\xi=\rho\,e_\xi$. Then

(a) the function $\rho$ has no strong maximums;

(b) there is a bivector field $e_0\wedge\nabla_{e_0}e_0$  such that $e_\xi$ is parallel along it.

Particulary, if $n=2$ then either $M^2$ is flat or $e_0$ is a geodesic vector field and $\rho$ is linear with
respect to the natural parameter along each $e_0$ geodesic line. Moreover, the field $\xi$ makes a constant
angle with each $e_0$ geodesic line.
\end{corollary}
\begin{proof}
Indeed, for any unit vector field $\eta$ consider the linear mapping $\nabla_Z\eta|_q : T_qM^n\to \eta^\perp_q$,
where $\eta^\perp_q$ is an orthogonal complement to $\eta$ in $T_qM^n$. For dimensional reasons it follows that
the kernel of this mapping is not empty. In other words, there exists a (unit) vector field $e_0$ such that
$\nabla_{e_0}\eta=0$.

Let $e_0$ be a unit vector field such that $\nabla_{e_0}e_\xi=0$. Then from $(\ref{Eq2})_1$ we conclude
$$
Hess_\rho(e_0,e_0)=0
$$
at each point of $M^n$. Therefore, the Hessian of $\rho$ can not be positively definite.

Moreover, from $(\ref{Eq2})_2$ we see that $r(e_0,e_0)e_\xi=0$, which gives
$\nabla_{e_0}\nabla_{e_0}e_\xi-\nabla_{\nabla_{e_0}e_0}e_\xi=-\nabla_{\nabla_{e_0}e_0}e_\xi=0$. Setting
$Z=e_0\wedge\nabla_{e_0}e_0$, we get
$
\nabla_Ze_\xi=0.
$

Suppose now that $n=2$. If $Z\ne0$ then $e_\xi$ is a parallel vector field on $M^2$ which means that $M^2$ is
flat. If $Z=0$ then evidently $e_0$ is a geodesic vector field. Since in this case $Hess_\rho
(e_0,e_0)=e_0(e_0(\rho))=0$,  we conclude that $\rho$ is linear with respect to the natural parameter along each
$e_0$ geodesic line.

As concerns the angle function $\big<e_0,e_\xi\big>$, we have
$$
e_0\big<e_0,e_\xi\big>=\big<\nabla_{e_0}e_0,e_\xi\big>+\big<e_0,\nabla_{e_0}e_\xi\big>=0.
$$
\end{proof}

Taking into account the Corollary \ref{Cor1}, introduce on $M^2$ a semi-geodesic coordinate system $(u,v)$ such
that $e_\xi$ is parallel along $u$-geodesics. Let
\begin{equation}\label{metric}
ds^2=du^2+b^2(u,v)\,dv^2
\end{equation}
be the first fundamental form of $M^2$ with respect to this coordinate system. Denote by $\partial_1$ and
$\partial_2$ the corresponding coordinate vector fields. Then the following equations should be satisfied:
$$
\nabla_{\partial_1}e_\xi=0, \qquad \partial_1^2(\rho)=0.
$$
Introduce the unit vector fields
$$
e_1=\partial_1,\qquad e_2=\frac{1}{b}\partial_2.
$$
Then the following rules of covariant derivation are valid
\begin{equation}\label{Frenet}
\begin{array}{ll}
\nabla_{e_1}e_1=0,\quad &\nabla_{e_1}e_2=0,\\[1ex]
\nabla_{e_2}e_1=-k\,e_2 \quad &\nabla_{e_2}e_2=k\,e_1,
\end{array}
\end{equation}
where $k$ is a (signed) geodesic curvature of $v$-curves. Remark that
$$
k=-\frac{\partial_1 b}{b}.
$$

With respect to chosen coordinate system, the field $\xi$ can be expressed as
\begin{equation}\label{field}
\xi=\rho\,(\cos \omega\,e_1+\sin\omega\,e_2),
\end{equation}
where $\omega=\omega(u,v)$ is an angle function, i.e.
$$
e_\xi=\cos \omega\,e_1+\sin\omega\,e_2.
$$
Introduce a unit vector field $\nu_\xi$ by
$$
\nu_\xi=-\sin \omega\,e_1+\cos\omega\,e_2.
$$
Then we can easily find
$$
\begin{array}{l}
\nabla_{e_1}e_\xi=\partial_1\omega\,\nu_\xi, \\[1ex]
\nabla_{e_2}e_\xi=(e_2(\omega)-k)\,\nu_\xi.
\end{array}
$$
Since $e_\xi$ is parallel along $u$-curves, we conclude that $\partial_1\omega=0$, so that $\omega=\omega(v)$.

Now the problem can be formulated as

{\it On a Riemannian 2-manifold with the metric (\ref{metric}), find a vector field of the form (\ref{field})
with
\begin{equation}\label{cond}
\partial^2_1\rho=0 \mbox{  and  } \omega=\omega(v)
\end{equation}
satisfying the equation (\ref{Eq3}).}

\begin{lemma}\label{tgcond1}
Let $M^2$ be a Riemannian 2-manifold with the metric (\ref{metric}) and $\xi$ be a local vector field on $M^2$
satisfying (\ref{cond}). Then $\xi$ is totally geodesic if and only if
\begin{equation}\label{tgEqn1}
\begin{array}{ll}
\nabla_{e_2}\nabla_{e_2}\xi-(k+c \,K)\nabla_{e_1}\xi=0,\\[1ex]
\nabla_{e_1}\nabla_{e_2}\xi+\nabla_{e_2}\nabla_{e_1}\xi+(k+c\, K)\nabla_{e_2}\xi=0,
\end{array}
\end{equation}
or in a scalar form
\begin{equation}\label{tgEqn2}
\begin{array}{l}
    \left\{ \begin{array}{l}
            e_2(e_2(\rho))-(k+c\,K)\,e_1(\rho)=\rho\lambda^2,\\[1ex]
            e_2(c)=0,
            \end{array}
    \right.\\[2ex]
    \left\{ \begin{array}{l}
            2e_1(e_2(\rho))+c\,K\,e_2(\rho)=0,\\[1ex]
            e_1(c)+c\,(k+c\,K)=0
            \end{array}
    \right.
\end{array}
\end{equation}
where $\lambda:=\big<\nabla_{e_2}e_\xi,\nu_\xi\big>=e_2(\omega)-k$, $c:=\rho^2\lambda=\pm
|\,\xi\wedge\nabla_{e_2}\xi|$ and $K$ is the Gaussian curvature of $M^2$.
\end{lemma}

\begin{proof}
Indeed,
$$
\begin{array}{l}
\nabla_{e_1}\xi=e_1(\rho)\,e_\xi,\\[1ex]
\nabla_{e_2}\xi=e_2(\rho)\,e_\xi+\rho\lambda\nu_\xi.
\end{array}
$$
So, taking into account (\ref{Frenet}) and (\ref{cond}), we have
$$
\begin{array}{l}
r(e_1,e_1)\xi=\nabla_{e_1}\nabla_{e_1}\xi-\nabla_{\nabla_{e_1}e_1}\xi=e_1(e_1(\rho))\,e_\xi=
\partial^2_1\rho\,e_\xi=0,\\[1ex]
r(e_1,e_2)\xi=\nabla_{e_1}\nabla_{e_2}\xi-\nabla_{\nabla_{e_1}e_2}\xi=\nabla_{e_1}\nabla_{e_2}\xi,\\[1ex]
r(e_2,e_1)\xi=\nabla_{e_2}\nabla_{e_1}\xi-\nabla_{\nabla_{e_2}e_1}\xi=\nabla_{e_2}\nabla_{e_1}\xi+
k\nabla_{e_2}\xi,\\[1ex]
r(e_2,e_2)\xi=\nabla_{e_2}\nabla_{e_2}\xi-\nabla_{\nabla_{e_2}e_2}\xi=\nabla_{e_2}\nabla_{e_2}\xi-
k\nabla_{e_1}\xi.
\end{array}
$$
As concerns the right-hand side of (\ref{Eq3}), we have
$$
\begin{array}{l}
R(\xi,\nabla_{e_1}\xi)e_1=0,\quad
R(\xi,\nabla_{e_1}\xi)e_2=0,\\[1ex]
R(\xi,\nabla_{e_2}\xi)e_1=\rho^2\lambda\,R(e_\xi,\nu_\xi)e_1=-\rho^2\lambda\,K\,e_2,\\[1ex]
R(\xi,\nabla_{e_2}\xi)e_2=\rho^2\lambda\,R(e_\xi,\nu_\xi)e_2=\rho^2\lambda\,K\,e_1.
\end{array}
$$
Therefore, setting $X=e_1$ in (\ref{Eq3}), we obtain an identity. Setting $X=e_2$, we have $$
 \nabla_{e_2}\nabla_{e_2}\xi-k\nabla_{e_1}\xi=\rho^2\lambda\,K\nabla_{e_1}\xi.
 $$
Setting $X=e_1+e_2$, we obtain
$$
r(e_1,e_2)\xi+r(e_2,e_1)\xi=-\rho^2\lambda\,K\,\nabla_{e_2}\xi,
$$
which can be reduced to
$$
\nabla_{e_1}\nabla_{e_2}\xi+\nabla_{e_2}\nabla_{e_1}\xi+k\nabla_{e_2}\xi=-\rho^2\lambda\, K\nabla_{e_2}\xi.
$$
It remains to mention that
$$
|\,\xi\wedge\nabla_{e_2}\xi|^2=|\xi|^2\,|\nabla_{e_2}\xi|^2-\big<\xi,\nabla_{e_2}\xi\big>^2=
\rho^2(e_2(\rho)^2+\rho^2\lambda^2)-(e_2(\rho)\rho)^2=\rho^4\lambda^2.
$$
So, if we set $c=\rho^2\lambda$, we evidently obtain (\ref{tgEqn1}).

Moreover, continuing calculations, we see that
$$
\begin{array}{rl}
\nabla_{e_2}\nabla_{e_2}\xi=\!\!&\Big[e_2(e_2(\rho))-\rho\lambda^2\Big]\,e_\xi+
\Big[e_2(\rho)\lambda+e_2(\rho\lambda)\Big]\,\nu_\xi=\\[1ex]
&\Big[e_2(e_2(\rho))-\rho\lambda^2\Big]\,e_\xi+\frac{1}{\rho}e_2(c)\,\nu_\xi,\\[2ex]
\nabla_{e_1}\nabla_{e_2}\xi+\nabla_{e_2}\nabla_{e_1}\xi=\!\!
&\!\!\Big[e_2(e_1(\rho))+e_1(e_2(\rho))\Big]e_\xi+\Big[e_1(\rho)\lambda+e_1(\rho\lambda)\Big]\nu_\xi=\\[1ex]
&\!\!\Big[e_2(e_1(\rho))+e_1(e_2(\rho))\Big]e_\xi+\frac{1}{\rho}e_1(c)\nu_\xi.
\end{array}
$$
Taking into account that $e_1(e_2(\rho))-e_2(e_1(\rho))=k\,e_2(\rho)$, the equations (\ref{tgEqn1}) can be
written as
$$
\begin{array}{l}
\Big[e_2(e_2(\rho))-\rho\lambda^2\Big]\,e_\xi+\frac{1}{\rho}e_2(c)\,\nu_\xi-(k+cK)e_1(\rho)\,e_\xi=0\\[1ex]
\Big[2\,e_1(e_2(\rho))-k\,e_2(\rho)\Big]\,e_\xi+\frac{1}{\rho}e_1(c)\,\nu_\xi+(k+cK)\Big[e_2(\rho)\,e_\xi+
\rho\lambda\,\nu_\xi\Big]=0
\end{array}
$$
and after evident simplifications we obtain the equations (\ref{tgEqn2}).
\end{proof}

\begin{proposition}
Let $M^2$ be a Riemannian manifold of constant curvature. Suppose $\xi$ is a non-zero local vector field on
$M^2$ such that $\xi(M^2)$ is totally geodesic in $TM^2$. Then $M^2$ is flat.
\end{proposition}
\begin{proof}
Let $M^2$ a Riemannian manifold of constant curvature $K\ne0$. Then the function $b$ in (\ref{metric}) should
satisfy the equation
$$
-\frac{\partial_{11}b}{b}=K.
$$
The general solution of this equation can be expressed in 3 forms:
\begin{itemize}
\item[(a)] $b(u,v)= A(v)\cos(u/r+\theta(v))$ or $b(u,v)=A(v)\sin(u/r+\theta(v))$ for $K=1/r^2>0$; \item[(b)]
$b(u,v)= A(v)\cosh(u/r+\theta(v))$ or $b(u,v)=A(v)\sinh(u/r+\theta(v))$ for $K=-1/r^2<0$; \item[(c)] $b(u,v)=
A(v)e^{u/r}$ for $K=-1/r^2<0$;
\end{itemize}

Evidently, we may set $A(v)\equiv 1$ (making a $v$-parameter change) in each of these cases.

 The equation $(\ref{tgEqn2})_2$ means that $c$ does not depend on $v$. Since $K$ is constant,
the equation $(\ref{tgEqn2})_4$ implies
$$
e_2(k)=0.
$$
If we remark that $k=-\frac{\partial_{1}b}{b}$ then one can easily find $\theta(v)=const$ in cases $(a)$ and
$(b)$.

 After a $u$-parameter change, the function $b$ takes one of the forms
\begin{itemize}
\item[(a)] $b(u,v)= \cos(u/r)$ or $b(u,v)=\sin(u/r)$ for $K=1/r^2>0$; \item[(b)] $b(u,v)= \cosh(u/r)$ or
$b(u,v)=\sinh(u/r)$ for $K=-1/r^2<0$; \item[(c)] $b(u,v)= e^{u/r}$ for $K=-1/r^2<0$;
\end{itemize}

From the equation $(\ref{tgEqn2})_4$ we find
$$
cK=-\frac{e_1(c)}{c}-k=-\frac{e_1(c)}{c}+\frac{e_1(b)}{b}=e_1(\ln b/c).
$$

Suppose first that $ e_2(\rho)\ne0$. Multiplying $(\ref{tgEqn2})_3$ by $e_2(\rho)$ we can easily solve this
equation with respect to $e_2(\rho)$ by a chain of simple transformations:
$$
\begin{array}{l}
2e_2(\rho)\cdot e_1(e_2(\rho))+e_1(\ln b/c)\cdot [e_2(\rho)^2]=0,\\[1ex]
e_1[e_2(\rho)^2]+e_1(\ln b/c)\cdot[e_2(\rho)^2]=0,\\[1ex]
\frac{e_1[e_2(\rho)^2]}{e_2(\rho)^2}+e_1(\ln b/c)=0, \\[1ex]
e_1[\ln e_2(\rho)^2]+e_1(\ln b/c)=0, \\[1ex]
 e_1( \ln[e_2(\rho)^2\,b/c])=0
\end{array}
$$
and therefore, $e_2(\rho)^2\,b/c=h(v)^2$ or
$$
 \partial_2\rho=h(v)\sqrt{c\,b}.
$$

Since $\rho$ is linear with respect to the $u$-parameter, say $\rho=a_1(v)u+a_2(v)$, then
$\partial_2\rho=a_1'u+a_2'$ and therefore $\sqrt{cb}$ is also linear with respect to $u$, namely
$\sqrt{cb}=m_1(v)u+m_2(v)=\frac{a_1'}{h}\,u+\frac{a_2'}{h}$. But the functions $c$ and $b$ do not depend on $v$.
Therefore $m_1$ and $m_2$ are constants, so $a_1=m_1\int h(v)\,dv, a_2=m_2\int h(v)\, dv$. Thus
$$
\sqrt{cb}=m_1u+m_2.
$$

Now the function $c$ takes the form
$$
c(u)=\frac{(m_1u+m_2)^2}{b}
$$
and therefore
$$
e_1(c)=\frac{2m_1(m_1u+m_2)}{b}-\frac{(m_1u+m_2)^2\partial_1 b}{b^2}.
$$
Substitution into $(\ref{tgEqn2})_4$ gives
$$
\frac{2m_1(m_1u+m_2)}{b}-\frac{2(m_1u+m_2)^2\partial_1 b}{b^2}+\frac{(m_1u+m_2)^4}{b^2}K=0
$$
or
$$
\frac{(m_1u+m_2)}{b^2}\Big[2m_1b-2(m_1u+m_2)\partial_1b+(m_1u+m_2)^3K \Big]=0.
$$
The expression in brackets is an algebraic one and can not be identically zero if $K\ne 0$. Therefore
$m_1=m_2=0$ and hence $\rho^2\lambda:= c=0$. But this identity implies $\lambda=0$ or $ \rho=0$. If $\lambda=0$
then $e_\xi$ is a parallel unit vector field and therefore, $M^2$ is flat and we come to a contradiction.
Therefore $\rho=0$.

\begin{remark}
If $K= 0$, we can not conclude that $c=0$. In this case the expression in brackets can be identically zero for
$m_1=0$ and $b=const$. And we have $c=m_2=const$.
\end{remark}

Suppose now that $e_2(\rho)=0$. Then
$$
\rho=a_1u+a_2,
$$
where $a_1,a_2$ are constants and we obtain the following system
\begin{equation}\label{tgEqn3}
\begin{array}{l}
-(k+cK)\partial_1\rho=\rho\lambda^2,\\[1ex]
\partial_2\,c=0, \\[1ex]
\partial_1c+c(k+cK)=0.
\end{array}
\end{equation}

If $\partial_1\rho=0$ then immediately $\rho=0$ or $\lambda=0$. The identity $\lambda=0$ implies $K=0$ as above.
Therefore, $\rho=0$.

Suppose $\partial_1\rho\ne 0$ or equivalently $a_1\ne 0$. Then from $(\ref{tgEqn3})_1$ we get
\begin{equation}\label{subs}
(k+cK)=-\frac{\rho\lambda^2}{a_1}
\end{equation}

 Since
$c=\rho^2\lambda$, from $(\ref{tgEqn3})_2$ we see that $\partial_2\lambda=0$ or
$\partial_2\left[\frac{\partial_2\omega+\partial_1b}{b}\right]=0$. Since $b$ does not depend on $v$, we have
$\partial_{22}\omega=0$ or equivalently $\partial_2\omega=\alpha=const$. Thus,
$\lambda=\frac{\alpha+\partial_1b}{b}$.

Now we can find $\partial_1c$ in two ways. First, from $(\ref{tgEqn3})_3$ using (\ref{subs}) and keeping in mind
that $c=\rho^2\lambda$:
$$
\partial_1c=c\frac{\rho\lambda^2}{a_1}=\frac{\rho^3\lambda^3}{a_1}
$$

Second, directly:
$$
\partial_1c=2\rho\partial_1\rho\lambda+\rho^2\partial_1\lambda.
$$
It is easy to see that $\partial_1\lambda=k\lambda-K$ and hence we get
$$
\partial_1c=2a_1\rho\lambda+\rho^2(k\lambda-K).
$$
Equalizing, we have
$$
2a_1\rho\lambda+\rho^2(k\lambda-K)-\frac{\rho^3\lambda^3}{a_1}=0
$$
or
$$
\frac{\rho}{a_1}\left[ 2a_1^2\lambda+a_1\rho(k\lambda-K)-\rho^2\lambda^3\right]=0.
$$
The expression in brackets is an algebraic one and can not be identically zero for $K\ne 0$. Since $\rho\not=
0$, we obtain a contradiction.

\begin{remark}
We do not obtain a contradiction if $K= 0$, since we have another solution $\lambda=0$ which gives
$\partial_1b+\alpha=0$ and hence  $b=-\alpha u+m$.
\end{remark}
\end{proof}

We have achieved the  result by putting a restriction on the geometry of the base manifold. Putting a
restriction on the vector field we are able to achieve a similar result. Recall that a totally geodesic vector
field necessarily makes a constant angle with some family of geodesics on the base manifold ( see Corollary
\ref{Cor1}). It is not parallel along this family and this fact is essential for its totally geodesic property.
Namely,

\begin{proposition}\label{Parallel}
Let $M^2$ be a Riemannian manifold. Suppose $\xi$ is a non-zero local vector field on $M^2$ which is parallel
along some family of geodesics of $M^2$. If $\xi(M^2)$ is totally geodesic in $TM^2$ then $M^2$ is flat.
\end{proposition}

\begin{remark}
Geometrically, this assertion means that if $\xi(M^2)$ is not transversal to the horizontal distribution on
$TM^2$  then $\xi(M^2)$ is never totally geodesic in $TM^2$ except when $M^2$ is flat.
\end{remark}

\begin{proof} Let $M^2$ be a \emph{non-flat} Riemannian manifold and suppose that the hypothesis of
the theorem is fulfilled. Then, choosing a coordinate system as in Lemma \ref{tgcond1}, we have
$$
\nabla_{e_1}\xi=0
$$
and we can reduce (\ref{tgEqn1}) to
\begin{equation}\label{parEqn}
\begin{array}{l}
\nabla_{e_2}\nabla_{e_2}\xi=0,\\
\nabla_{e_1}\nabla_{e_2}\xi+(k+cK)\nabla_{e_2}\xi=0.
\end{array}
\end{equation}

Now  make a simple computation.
$$
\begin{array}{l}
R(e_2,e_1)\nabla_{e_2}\xi=\nabla_{e_2}\nabla_{e_1}\nabla_{e_2}\xi-
\nabla_{e_1}\nabla_{e_2}\nabla_{e_2}\xi-\nabla_{[e_2,e_1]}\nabla_{e_2}\xi=\\[1ex]
\nabla_{e_2}\nabla_{e_1}\nabla_{e_2}\xi-k\nabla_{e_2}\nabla_{e_2}\xi= \nabla_{e_2}\nabla_{e_1}\nabla_{e_2}\xi.
\end{array}
$$
On the other hand, differentiating $(\ref{parEqn})_2$, we find
$$
\nabla_{e_2}\nabla_{e_1}\nabla_{e_2}\xi=-e_2(k+cK)\nabla_{e_2}\xi.
$$
So we have
$$
R(e_2,e_1)\nabla_{e_2}\xi=-e_2(k+cK)\nabla_{e_2}\xi.
$$
Therefore, either $\nabla_{e_2}\xi=0$ or $e_2(k+cK)=0$. If we accept the first case we see that $\xi$ is a
parallel vector field on $M^2$ and we get a contradiction.

If we accept the second case, we obtain
$$
R(e_2,e_1)\nabla_{e_2}\xi=0,
$$
which means that $\nabla_{e_2}\xi$ belongs to a kernel of the curvature operator of $M^2$. In dimension 2 this
means that $M^2$ is flat or, equivalently, $\xi$ is a parallel vector field and we obtain a contradiction, as
well.

\end{proof}

\subsection{Ruled totally geodesic submanifolds in $TM^2$}\label{Ruled}

\begin{proposition}
Let $M^2$ be a Riemannian manifold of sign-preserving curvature. Consider a ruled submanifold  $\tilde F^2$  in
$TM^2$ given locally  by
$$
\gamma(s): \left\{\begin{array}{l}
                x^1=x^1(s),\\[1ex]
                x^2=x^2(s),
                \end{array}
            \right. \qquad
 \xi(t,s): \left\{
           \begin{array}{l}
            \xi^1=t,\\[1ex]
            \xi^2=\xi^2(t,s).
            \end{array}
            \right.
$$
Then $\tilde F^2$ is totally geodesic in $TM^2$ if $\gamma(s)$ is a geodesic in $M^2$,
$$\xi(t,s)=t\,\rho(s)\,e(s),$$ where $e(s)$ is a unit vector field which is parallel along $\gamma$ and $\rho(s)$
is an arbitrary smooth function.
\end{proposition}
\begin{remark} Geometrically, $\tilde F^2$ is a cylinder-type surface based on geodesic $\gamma(s)$ with elements
directed by a unit vector field $e(s)$ parallel along $\gamma(s)$.
\end{remark}
\begin{proof}
Fixing $s=s_0$, we see that $F^2$ meets the fiber over $x^1(s_0),x^2(s_0)$ by a curve $\xi(t,s_0)$. If $F^2$ is
supposed to be totally geodesic, then this curve is a straight line on the fiber. Therefore, the family
$\xi(t,s)$ should be of the form
$$
\xi(t,s): \left\{
\begin{array}{l}
\xi^1=t,\\[1ex]
\xi^2=\alpha(s)t+\beta(s).
\end{array}
\right.
$$
Introduce two vector fields given along $\gamma(s)$ by
\begin{equation}\label{fields1}
a=\pt_1+\alpha(s)\pt_2,\quad b=\beta(s)\pt_2.
\end{equation}
Then we can represent $\xi(t,s)$ as
$$
\xi(t,s)=a(s)\,t+b(s).
$$

Denote by $\tau$ and $\nu$ the vectors of the Frenet frame of the curve $\gamma(s)$. Denote also by (\,$'$) the
covariant derivative of vector fields with respect to the arc-length parameter on $\gamma(s)$. Then
$$
\left\{
\begin{array}{l}
\tau\,'=k\,\nu,\\
\nu\,'=-k\,\tau.
\end{array}
\right.
$$

Denote by $\tilde\pt_1, \tilde\pt_2$ the $s$ and $t$ coordinate vector fields on $F^2$ respectively. A simple
calculation yields
$$
\tilde\pt_1={\tau}^h+(\xi')^v,\quad \tilde\pt_2=a^v.
$$

One of the unit normal vector fields can be found immediately, namely $\tilde N_1={\nu\,}^h$. Consider the
conditions on $F^2$ to be totally geodesic with respect to the normal vector field $\tilde N_1$. Using
formulas (\ref{Kow}),
$$
\tilde\nabla_{\tilde \pt_1}\tilde N_1=\tilde\nabla_{\tau\,^h+(\xi')^v}\nu\,^h=
-k\tau\,^h-\frac12\Big[R(\tau,\nu)\xi\Big]^v+\frac12\Big[R(\xi,\xi')\nu\Big]^h
$$
Therefore,
$$
\big<\big<\tilde\nabla_{\tilde \pt_1}\tilde N_1,\tilde\pt_2\big>\big>= -\frac12\big<R(\tau,\nu)\xi,a\big>=
-\frac12\big<R(\tau,\nu)b,a\big>=0.
$$
Since $M^2$ is supposed to be non-flat, it follows $b\wedge a=0$. From (\ref{fields1}) we conclude $b=0$. Thus,
$\xi(t,s)=a(s)\,t$. Moreover,
$$
\begin{array}{rl}
\big<\big<\tilde\nabla_{\tilde \pt_1}\tilde N_1,\tilde\pt_1\big>\big>=&-k-\frac12\big<R(\tau,\nu)\xi,\xi'\big>
+\frac12\big<R(\xi,\xi')\nu,\tau\big>=\\[1ex]
&-k+\big<R(\xi,\xi')\nu,\tau\big>=-k+t^2\big<R(a,a')\nu,\tau\big>=0
\end{array}
$$
identically with respect to parameter $t$. Therefore, $k=0$ and $a\wedge a'=0$. Thus, $\gamma(s)$ is a geodesic
line on $M^2$. In addition,  $(a\wedge a'=0)\sim (a'=\lambda a)$. Set $a=\rho(s)\, e(s)$, where $\rho=|a(s)|$.
Then $(a'=\lambda a)\sim(\rho'\,e+\rho\,e'=\lambda\rho\,e) $, which means that $e'=0$. From this we conclude
$$
\xi(t,s)=t\rho(s)\,e(s),
$$
where $\rho(s)$ is arbitrary function and $e(s)$ is a unit vector field, parallel along $\gamma(s)$. Therefore,
$$
\tilde\pt_1={\tau}^h+t\rho\,'\,e^v,\quad \tilde\pt_2=\rho \,e^v
$$
and we can find another unit normal vector field $\tilde N_2=(e^\perp)^v$, where $e^\perp(s)$ is a unit vector
field also parallel along $\gamma(s)$ and orthogonal to $e(s)$. For this vector field we have
$$
\begin{array}{l}
\tilde\nabla_{\tilde \pt_1}\tilde N_2=\tilde\nabla_{\tau\,^h+(\xi')^v}(e^\perp)^v=
[(e^\perp)']^v+\frac12\Big[R(\xi,e^\perp)\tau\Big]^h=\frac12t\rho \Big[R(e,e^\perp)\tau\Big]^h,\\[1ex]
\tilde\nabla_{\tilde \pt_2}\tilde N_2=0
\end{array}
$$
Evidently, $ \big<\big<\tilde\nabla_{\tilde \pt_i}\tilde N_2,\tilde\pt_k\big>\big>=0$ for all $i,k=1,2$. Thus,
the submanifold is totally geodesic.
\end{proof}
The converse statement is true in general.
\begin{proposition}\label{Gener}
Let $M^n$ be a Riemannian manifold. Consider a cylinder type surface $\tilde F^2\subset TM^n$  parameterized as
$$
\big\{\gamma(s),t\,\rho(s)\,e(s)\big\},
$$
where $\gamma(s)$ is a geodesic in $M^n$, $e(s)$ is a unit vector field, parallel along $\gamma$ and $\rho(s)$
is an arbitrary smooth function. Then $\tilde F^2$ is totally geodesic in $TM^n$ and intrinsically flat.
\end{proposition}
\begin{proof}
Indeed, the tangent basis of $\tilde F^2$ is consisted of
$$
\tilde \pt_1={\gamma\,'}^h+t\rho\,'e^v,\qquad \tilde \pt_2=\rho\, e^v.
$$
By   formulas (\ref{Kow}),
$$
\begin{array}{l}
\ds \tilde \nabla_{\tilde \pt_1}{\tilde \pt_1}=(\nabla_{\gamma\,'}{\gamma\,'})^h+
\frac12\big[R(\gamma\,',\gamma\,')\xi\big]^v=0\,,\\[2ex]
\ds \tilde \nabla_{\tilde \pt_1}{\tilde \pt_2}=(\nabla_{\gamma\,'}\rho\,e)^v+
\frac12\big[R(\xi,\rho\,e)\gamma\,'\big]^h=\rho\,'e^v\ \sim \ \tilde\pt_2\,,\\[2ex]
\ds \tilde \nabla_{\tilde \pt_2}{\tilde \pt_1}=\frac12\big[R(\xi,\rho\,e)\gamma\,'\big]^h=
\frac12\big[R(t\,\rho\,e ,\rho\,e)\gamma\,'\big]^h=0\\[2ex]
\ds \tilde \nabla_{\tilde \pt_2}{\tilde \pt_2}=\rho e^v(\rho)\,e^v=0.
\end{array}
$$
It is easy to find the Gaussian curvature of this submanifold, since it is equal to the sectional curvature of
$TM^2$ along the $\tilde\pt_1\wedge\tilde\pt_2$- plane. Using the curvature tensor expressions \cite{Kow}, we
find
$$
Gauss(\tilde F^2)=\big<\big<\tilde R(\tau\,^h,e^v)e^v,\tau\,^h\big>\big>=\frac14|R(\xi,e)\tau|^2=0.
$$

\end{proof}

\section{Local description of 3-dimensional totally geodesic submanifolds in $TM^2$}\label{3-dim}

\begin{theorem}
Let $M^2$ be Riemannian manifold with Gaussian curvature $K$. A totally geodesic submanifold $\tilde F^3\subset
TM^2$ locally is either

a) a 3-plane in $TM^2=E^4$ if $K=0$, or

b) a restriction of the tangent bundle to a geodesic $\gamma\in M^2$ such that $K|_\gamma=0$ if $K\not\equiv0$.
If $M^2$ does not contain such a geodesic, then $TM^2$ does not admit 3-dimensional totally geodesic
submanifolds.
\end{theorem}

\begin{proof}
Let $\tilde F^3$ be a submanifold it $TM^2$. Let $(x^1,x^2;\xi^1,\xi^2)$ be a local chart on $TM^2$. Then
locally $\tilde F^3$ can be given mapping $f$ of the form
$$
f: \left\{
\begin{array}{l}
x^1=x^1(u^1,u^2,u^3),\\[1ex]
x^2=x^2(u^1,u^2,u^3),\\[1ex]
\xi^1=\xi^1(u^1,u^2,u^3),\\[1ex]
\xi^2=\xi^2(u^1,u^2,u^3),\\[1ex]
\end{array}
 \right.
$$
where $u^1,u^2,u^3$ are the local parameters on $\tilde F^3$. The Jacobian matrix $f_*$ of the mapping $f$ is of
the form {\large
$$
f_*=
\left(
\begin{array}{ccc}
 \frac{\pt x^1}{\pt u^1} & \frac{\pt x^1}{\pt u^2}&  \frac{\pt x^1}{\pt u^3} \\[1ex]
 \frac{\pt x^2}{\pt u^1} & \frac{\pt x^2}{\pt u^2}&  \frac{\pt x^2}{\pt u^3} \\[1ex]
\frac{\pt \xi^1}{\pt u^1} & \frac{\pt \xi^1}{\pt u^2} & \frac{\pt \xi^1}{\pt u^3}  \\[1ex]
 \frac{\pt \xi^2}{\pt u^1} & \frac{\pt \xi^2}{\pt u^2}&  \frac{\pt \xi^2}{\pt u^3} \\[1ex]
\end{array}
\right).
$$
} Since $rank \ f_*=3$, we have two geometrically different possibilities to achieve the rank, namely {\large
$$
(a)\quad
\det
\left(%
\begin{array}{ccc}
 \frac{\pt x^1}{\pt u^1} & \frac{\pt x^1}{\pt u^2}&  \frac{\pt x^1}{\pt u^3} \\[1ex]
 \frac{\pt x^2}{\pt u^1} & \frac{\pt x^2}{\pt u^2}&  \frac{\pt x^2}{\pt u^3} \\[1ex]
\frac{\pt \xi^1}{\pt u^1} & \frac{\pt \xi^1}{\pt u^2} & \frac{\pt \xi^1}{\pt u^3}  \\[1ex]
\end{array}
\right)\ne0; \quad  (b)\ \ \det
 \left(%
\begin{array}{ccc}
 \frac{\pt x^1}{\pt u^1} & \frac{\pt x^1}{\pt u^2}&  \frac{\pt x^1}{\pt u^3} \\[1ex]
\frac{\pt \xi^1}{\pt u^1} & \frac{\pt \xi^1}{\pt u^2} & \frac{\pt \xi^1}{\pt u^3}  \\[1ex]
 \frac{\pt \xi^2}{\pt u^1} & \frac{\pt \xi^2}{\pt u^2}&  \frac{\pt \xi^2}{\pt u^3} \\[1ex]
\end{array}%
\right)\ne0.
$$
} Without loss  of generality we can consider this possibilities in such a way that (b) excludes (a).

\textbf{Consider the case (a).}  In this case we can locally parameterize the submanifold $F^3$ as
$$
f:\
\left\{
\begin{array}{l}
x^1=u^1,\\[1ex]
x^2=u^2,\\[1ex]
\xi^1=u^3,\\[1ex]
\xi^2=\xi^2(u^1,u^2,u^3).
\end{array}
\right.
$$
By hypothesis, the submanifold $\tilde F^3$ is totally geodesic in $TM^2$. Therefore, it intersects each fiber
of $TM^2$ by a vertical geodesic, i.e. by a straight line. Fix $u_0=(u^1_0, \, u^2_0)$. Then the parametric
equation of $\tilde F^3\cap T_{u_0}M^2$  with respect to fiber parameters is
$$\left\{
\begin{array}{l}
\xi^1=u^3,\\[1ex]
\xi^2=\xi^2(u^1_0,u^2_0,u^3).
\end{array}
\right.
$$
On the other hand, this equation should be the equation of a straight line and hence
$$\left\{
\begin{array}{l}
\xi^1=u^3,\\[1ex]
\xi^2=\alpha(u^1_0,u^2_0)\,u^3+\beta(u^1_0,u^2_0),
\end{array}
\right.
$$
where $\alpha(u)=\alpha(u^1,u^2)$ and $\beta(u)=\beta(u^1,u^2)$ some smooth functions on $M^2$. From this
viewpoint, after setting $u^3=t$  \textit{the submanifold under consideration can be locally represented as a
one-parametric family of smooth vector fields $\xi_t$ on $M^2$} of the form
$$
\xi_t(u)=t\,\pt_1+\big(\alpha(u)t+\beta(u)\big)\,\pt_2
$$
with respect to the coordinate frame $\pt_1=\pt/\pt u^1,\,\pt_2=\pt/\pt u^2$.

Introduce the vector fields
\begin{equation}\label{fields2}
a(u)=\pt_1+\alpha(u)\,\pt_2, \quad b(u)=\beta(u)\,\pt_2.
\end{equation}
Then $\xi_t$ can be expressed as
$$
\xi_t(u)=t\,a(u)+b(u).
$$
It is natural to denote by $\xi_t(M^2)$ a submanifold $\tilde F^3\subset TM^2$ of this kind.

Denote by $\tilde \pt_i \ (i=1,\dots,3)$ the coordinate vector fields of $\xi_t(M^2)$. Then
$$
\begin{array}{l}
  \tilde \pt_1=\big\{1,0,0,t\,\pt_1\alpha+\pt_1\beta\big\},\\[1ex]
  \tilde \pt_2=\big\{0,1,0,t\,\pt_2\alpha+\pt_2\beta\big\},\\[1ex]
  \tilde \pt_3=\big\{0,0,1,\alpha\big\}.\\[1ex]
\end{array}
$$
A direct calculation shows that these fields can be represented as
$$
\begin{array}{l}
\tilde \pt_1=\pt_1^h+t(\nabla_{\pt_1}\,a)^v+(\nabla_{\pt_1}\,b)^v,\\[1ex]
\tilde \pt_2=\pt_2^h+t(\nabla_{\pt_2}\,a)^v+(\nabla_{\pt_2}\,b)^v,\\[1ex]
\tilde \pt_3=a^v.
\end{array}
$$
Denote by $\tilde N$ a normal vector field of $\xi_t(M^2)$. Then
$$
\tilde N=(a^\perp)^v+Z_t^h,
$$
where $\big<a^\perp,a\big>=0$ and the field $Z_t=Z_t^1\pt_1+Z_t^2\pt_2$ can be found easily from the equations
$$
\begin{array}{l}
\big<\big<\tilde \pt_i,\tilde
N\big>\big>=\big<Z_t,\pt_i\big>+t\big<\nabla_{\pt_i}\,a,a^\perp\big>+\big<\nabla_{\pt_i}\,b,a^\perp\big>=0
\quad (i=1,2)
\end{array}
$$
Using the   formulas (\ref{Kow}), one can find
$$
\begin{array}{rl}
\tilde\nabla_{\tilde\pt_i}a^v=&\tilde \nabla_{\pt_i^h+t(\nabla_{\pt_i}\,a)^v+(\nabla_{\pt_i}\,b)^v }a^v=\\[1ex]
&(\nabla_{\pt_i}a)^v+\frac12\Big[R(\xi_t,a)\pt_i\Big]^h= (\nabla_{\pt_i}a)^v+\frac12\Big[R(b,a)\pt_i\Big]^h.
\end{array}
$$
If the submanifold $\xi_t(M^2)$ is totally geodesic, then the following equations should be satisfied
identically
$$
\big<\big<\tilde\nabla_{\tilde\pt_i}\tilde\pt_3,\tilde N\big>\big>=\big<\nabla_{\pt_i}a,a^\perp\big>+
\frac12\big<R(b,a)\pt_i,Z_t\big>=0
$$
with respect to the parameter $t$. To simplify the further calculations, suppose that  the coordinate system on
$M^2$ is the orthogonal one, so that $\big<\pt_1,\pt_2\big>=0$ and
$$
R(b,a)\pt_2=g^{11}K\,|b\wedge a|\pt_1,\quad R(b,a)\pt_1=-g^{22}K\,|b\wedge a|\pt_2,
$$
where $K$ is the Gaussian curvature of $M^2$ and $g^{11}, g^{22}$ are the contravariant metric coefficients.
Then we have
$$
\begin{array}{l}
\big<R(b,a)\pt_1,Z_t\big>=-g^{22}K\,|b\wedge a|\big<Z_t,\pt_2\big>=\\[1ex]
\hphantom{\big<R(b,a)\pt_1,Z_t\big>=} g^{22}K\,|b\wedge a|\Big(t\big<\nabla_{\pt_2}\,a,a^\perp\big>+\big<\nabla_{\pt_2}\,b,a^\perp\big>\Big),\\[2ex]
\big<R(b,a)\pt_2,Z_t\big>=g^{11}K\,|b\wedge a|\big<Z_t,\pt_1\big>= \\[1ex]
\hphantom{\big<R(b,a)\pt_1,Z_t\big>=} -g^{11}K\,|b\wedge
a|\Big(t\big<\nabla_{\pt_1}\,a,a^\perp\big>+\big<\nabla_{\pt_1}\,b,a^\perp\big>\Big).
\end{array}
$$
 Thus we get the system
$$
\left\{\begin{array}{l} g^{22}K|b\wedge a|\big<\nabla_{\pt_2}a,a^\perp\big>t+\big<\nabla_{\pt_1}a,a^\perp\big>+
g^{22}K|b\wedge a|\big<\nabla_{\pt_2}b,a^\perp\big>=0,\\[2ex]
g^{11}K|b\wedge a|\big<\nabla_{\pt_1}a,a^\perp\big>t-\big<\nabla_{\pt_2}a,a^\perp\big>+ g^{11}K|b\wedge
a|\big<\nabla_{\pt_1}b,a^\perp\big>=0,
\end{array}
\right.
$$
which should be satisfied identically with respect to $t$. As a consequence, we have 3 cases:
\begin{itemize}
\item[\bf(i)]\quad  $K=0,\ \left\{\begin{array}{l}\big<\nabla_{\pt_1}\,a,a^\perp\big>=0, \\[1ex]
\big<\nabla_{\pt_2}\,a,a^\perp\big>=0\end{array}\right.;$

\item[\bf(ii)]\quad $K\ne0$,\quad $|b\wedge a|=0,\ \left\{\begin{array}{l}\big<\nabla_{\pt_1}\,a,a^\perp\big>=0, \\[1ex]
\big<\nabla_{\pt_2}\,a,a^\perp\big>=0\end{array}\right.;$

\item[\bf(iii)]\quad $K\ne0$,\ \ $|b\wedge a|\ne0$,\ \ $\left\{\begin{array}{l}\big<\nabla_{\pt_1}\,a,a^\perp\big>=0, \\[1ex]
\big<\nabla_{\pt_2}\,a,a^\perp\big>=0\end{array}\right.,
\left\{\begin{array}{l}\big<\nabla_{\pt_1}\,b,a^\perp\big>=0, \\[1ex]
\big<\nabla_{\pt_2}\,b,a^\perp\big>=0\end{array}\right.;$
\end{itemize}

{Case (i).} In this case the base manifold is flat and we can choose a Cartesian coordinate system, so that the
covariant derivation becomes a usual one and we have
$$
\left\{
\begin{array}{l}
\nabla_{\pt_i}\,a=\big\{0,\pt_i\alpha\big\} \quad(i=1,2)\\[1ex]
a^\perp=\big\{-\alpha,1\big\}
\end{array}
\right.
$$
From $\big<\nabla_{\pt_i}\,a,a^\perp\big>=0$  it follows that $\alpha=const$, i.e. $a$ is a parallel vector
field. Moreover, in this case
$$
\begin{array}{l}
  \tilde \pt_1=\big\{1,0,0,\pt_1\beta\big\}=\pt_1^h+(\pt_1\,b)^v,\\[1ex]
  \tilde \pt_2=\big\{0,1,0,\pt_2\beta\big\}=\pt_1^h+(\pt_1\,b)^v,\\[1ex]
  \tilde \pt_3=\big\{0,0,1,\alpha\big\},\\[1ex]
  \tilde N=\big\{-\pt_1\beta,-\pt_2\beta,-\alpha,\,1\big\}.
\end{array}
$$
Now we can find
$$
\tilde\nabla_{\tilde\pt_i}\tilde\pt_k=(\nabla_{\pt_i}\pt_k b)^v=\big\{0,0,0,\pt_{ik}\beta\big\}
$$
and the conditions
$$
\big<\big<\tilde\nabla_{\tilde\pt_i}\tilde\pt_k,\tilde N\big>\big>=0
$$
imply $\pt_{ik}\beta=0$. Thus, $\beta=m_1u^1+m_2u^2+m_0$, where $m_1, m_2, m_0$ are arbitrary constants. As a
consequence, the submanifold $\xi_t(M^2)$ is described by parametric equations of the form
$$
\left\{
\begin{array}{l}
x^1=u^1,\\[1ex]
x^2=u^2,\\[1ex]
\xi^1=t,\\[1ex]
\xi^2=\alpha t+m_1u^1+m_2u^2+m_0
\end{array}
\right.
$$
and we have a hyperplane in $TM^2=E^4$. \vspace{2ex}

{Case(ii).} Keeping in mind (\ref{fields2}), the condition $b\wedge a=0$ implies $b=0$. The conditions
$$
\left\{\begin{array}{l}\big<\nabla_{\pt_1}\,a,a^\perp\big>=0, \\[1ex]
\big<\nabla_{\pt_2}\,a,a^\perp\big>=0\end{array}\right.
$$
imply $\nabla_{\pt_1}\,a=\lambda_1(u)\,a,\ \nabla_{\pt_2}\,a=\lambda_2(u)\,a$. As a consequence, we have
$$
\begin{array}{l}
\xi_t=t\,a\\[1ex]
  \tilde \pt_1=\pt_1^h+t(\nabla_{\pt_1}\,a)^v=\pt_1^h+t\lambda_1\,a^v,\\[1ex]
  \tilde \pt_2=\pt_2^h+t(\nabla_{\pt_2}\,a)^v=\pt_2^h+t\lambda_2\,a^v,\\[1ex]
  \tilde \pt_3=a^v,\\[1ex]
  \tilde N=(a^\perp)^v.
\end{array}
$$
Using   formulas (\ref{Kow}),
$$
\begin{array}{ll}
\tilde\nabla_{\tilde\pt_i}\tilde\pt_k=&\tilde\nabla_{\pt_i^h+t\lambda_i\,a^v}\Big(\pt_k^h+t\lambda_k\,a^v\Big)=\\[2ex]
&\tilde\nabla_{\pt_i^h}\pt_k^h+t\lambda_i\tilde\nabla_{a^v}\pt_k^h+\tilde\nabla_{\pt_i^h}(t\lambda_ka^v)+
t^2\lambda_i\lambda_k\tilde\nabla_{a^v}a^v=\\[2ex]
&(\nabla_{\pt_i}\pt_k)^h-\frac12\Big[R(\pt_i,\pt_k)\xi_t\Big]^v+t\lambda_i\frac12\Big[R(\xi_t,a)\pt_k\Big]^h+\\[1ex]
&t\pt_i(\lambda_k)a^v+t\lambda_k(\nabla_{\pt_i}a)^v+t\lambda_k\frac12\Big[R(\xi_t,a)\pt_i\Big]^h =\\[2ex]
&(\nabla_{\pt_i}\pt_k)^h-t\frac12\Big[R(\pt_i,\pt_k)a\Big]^v+t\pt_i(\lambda_k)a^v+t\lambda_k\lambda_i\,a^v.
\end{array}
$$
Evidently, for $i\ne k$
$$
\big<\big<\tilde\nabla_{\tilde\pt_i}\tilde\pt_k,\tilde
N\big>\big>=-t\frac12\big<R(\pt_i,\pt_k)a,a^\perp\big>\ne0,
$$
since $M^2$ is non-flat and $a\ne0$. Contradiction.\vspace{2ex}

{Case (iii).} The conditions imply
$$
\nabla_i a=\lambda_i(u)\, a,\quad \nabla_i b=\mu_i(u)\, a \quad (i=1,2)
$$
and we have
$$
\begin{array}{l}
\xi_t=t\,a+b\\[1ex]
  \tilde \pt_1=\pt_1^h+(t\lambda_1+\mu_1)\,a^v,\\[1ex]
  \tilde \pt_2=\pt_1^h+(t\lambda_2+\mu_2)\,a^v,\\[1ex]
  \tilde \pt_3=a^v,\\[1ex]
  \tilde N=(a^\perp)^v.
\end{array}
$$
A calculation as above leads to the identity
\begin{multline*}
$$
\ds \big<\big<\tilde\nabla_{\tilde\pt_i}\tilde\pt_k,\tilde
N\big>\big>=-\frac12\big<R(\pt_i,\pt_k)\xi_t,a^\perp\big>=\\
-t\frac12\big<R(\pt_i,\pt_k)a,a^\perp\big>-\frac12\big<R(\pt_i,\pt_k)b,a^\perp\big>=0
$$
\end{multline*}
which can be true if and only if
$$
\left\{ \begin{array}{l}
 \big<R(\pt_i,\pt_k)a,a^\perp\big>=0,\\[1ex]
 \big<R(\pt_i,\pt_k)b,a^\perp\big>=0.
 \end{array}
 \right.
$$
The first condition contradicts $K\ne0$.

\textbf{Consider the case (b).} In this case the submanifold $\tilde F^3$ can be locally parametrized  by
$$
\left\{
\begin{array}{l}
x^1=u^1,\\[1ex]
x^2=x^2(u^1,u^2,u^3)\\[1ex]
\xi^1=u^2,\\[1ex]
\xi^2=u^3.
\end{array}
\right.
$$
Since we exclude the case (a), we should suppose
$$
\det
\left(%
\begin{array}{ccc}
 \frac{\pt x^1}{\pt u^1} & \frac{\pt x^1}{\pt u^2}&  \frac{\pt x^1}{\pt u^3} \\[1ex]
 \frac{\pt x^2}{\pt u^1} & \frac{\pt x^2}{\pt u^2}&  \frac{\pt x^2}{\pt u^3} \\[1ex]
\frac{\pt \xi^1}{\pt u^1} & \frac{\pt \xi^1}{\pt u^2} & \frac{\pt \xi^1}{\pt u^3}  \\[1ex]
\end{array}
\right)=
\det
 \left(%
\begin{array}{ccc}
1 & 0&  0 \\[1ex]
 \frac{\pt x^2}{\pt u^1} & \frac{\pt x^2}{\pt u^2}&  \frac{\pt x^2}{\pt u^3} \\[1ex]
        0 & 1 & 0  \\[1ex]
\end{array}
\right)=-\frac{\pt x^2}{\pt u^3}=0;
$$
$$
\det
 \left(%
\begin{array}{ccc}
 \frac{\pt x^1}{\pt u^1} & \frac{\pt x^1}{\pt u^2}&  \frac{\pt x^1}{\pt u^3} \\[1ex]
 \frac{\pt x^2}{\pt u^1} & \frac{\pt x^2}{\pt u^2}&  \frac{\pt x^2}{\pt u^3} \\[1ex]
\frac{\pt \xi^2}{\pt u^1} & \frac{\pt \xi^2}{\pt u^2} & \frac{\pt \xi^2}{\pt u^3}  \\[1ex]
\end{array}
\right)= \det
 \left(%
\begin{array}{ccc}
1 & 0&  0 \\[1ex]
 \frac{\pt x^2}{\pt u^1} & \frac{\pt x^2}{\pt u^2}&  \frac{\pt x^2}{\pt u^3} \\[1ex]
        0 & 0 & 1  \\[1ex]
\end{array}
\right)=\frac{\pt x^2}{\pt u^2}=0;
$$
Therefore, in this case we have a submanifold, which can be parametrized by
$$
\left\{
\begin{array}{l}
x^1=x^1(s),\\[1ex]
x^2=x^2(s)\\[1ex]
\xi^1=u^2,\\[1ex]
\xi^2=u^3,
\end{array}
\right.
$$
where $s$ is a natural parameter of the  regular curve $\gamma(s)=\big\{x^1(s), x^2(s)\big\}$ on $M^2$.
Geometrically, a submanifold of this class is nothing else but the restriction of $TM^2$ to the curve
$\gamma(s)$. Denote by $\tau$ and $\nu$ the Frenet frame of $\gamma(s)$. It is easy to verify that
$$
 \tilde \pt_1=\tau\,^h,\quad  \tilde \pt_2=\pt_1^v,\quad  \tilde \pt_3=\pt_2^v,\quad
 \tilde N=\nu\,^h.
$$
By   formulas (\ref{Kow}), for $i=1,2$
$$
\big<\big<\tilde \nabla_{\tilde\pt_{1+i}}\tilde N,\tilde \pt_1 \big>\big>= \big<\big<\tilde
\nabla_{\pt_{i}^v}\nu^h,\tau^h
\big>\big>=\frac12\big<R(\xi,\pt_i)\nu,\tau\big>=\frac12\big<R(\tau,\nu)\pt_i,\xi\big>=0
$$
for arbitrary $\xi$. Evidently, $M^2$ must be flat along $\gamma(s)$.
\end{proof}

\end{document}